\documentclass{gtart}

\def\ifplaintex{\expandafter\ifx\csname documentclass\endcsname\relax}

\ifplaintex 
\else
\fi

\def\gtm{{\mathsurround=0pt\it $\cal G\mskip-2mu$eometry \&\ 
$\cal T\!\!$opology $\cal M\mskip-1mu$onographs}}    

\def\gtp{{\mathsurround=0pt\it $\cal G\mskip-2mu$eometry \&\ 
$\cal T\!\!$opology $\cal P\!$ublications}}  

\def\recd{{\small Received:\qua\receiveddate\ifx\reviseddate\relax
\else\qquad Revised:\qua\reviseddate\fi\par}} 

\def\volumenumber#1{\def\thevolumenumber{#1}}
\def\volumeyear#1{\def\thevolumeyear{#1}}
\def\volumename#1{\def\thevolumename{#1}}
\def\papernumber#1{\def\thepapernumber{#1}}
\def\pagenumbers#1#2{\def\startpage{#1}\def\finishpage{#2}}
\def\published#1{\def\publishdate{#1}}
\def\received#1{\def\receiveddate{#1}}
\def\revised#1{\def\reviseddate{#1}}
\def\accepted#1{\def\accepteddate{#1}}

\def\coverauthors#1{\def\thecoverauthors{#1}}
\def\asciiauthors#1{\def\theasciiauthors{#1}}

\def\coverauthors#1{\def\thecoverauthors{#1}}
\long\def\asciiabstract#1{\long\def\theasciiabstract{#1}}

\let\\\par
\let\thevolumenumber\relax\let\thepapernumber\relax
\let\thevolumeyear\relax\let\startpage\relax
\let\finishpage\relax\let\publishdate\relax\let\receiveddate\relax
\let\reviseddate\relax\let\accepteddate\relax\let\theasciititle\relax
\let\theasciiauthors\relax
\let\theasciiabstract\relax
\let\thecoverauthors\relax
\let\thecoverauthors\relax\let\theerratum\relax\let\theasciiemail\relax
\let\theshortauthors\relax\let\theshorttitle\relax

\def\startpage{1}\def\finishpage{15}\def\thepapernumber{77}

\volumenumber{2}
\volumename{Proceedings of the Kirbyfest}
\volumeyear{1999}

\long\def\maketitlep{   

\count0=\startpage

\gtm\nl        
{\small Volume \thevolumenumber: \thevolumename\nl 
\ifx\theerratum\relax\else Erratum \erratumnumber\nl\fi
Pages \startpage--\finishpage\nl}

\vglue 0.1truein   

{\parskip=0pt\leftskip 0pt plus 1fil\def\\{\par\smallskip}{\ifplaintex\large
\else\Large\fi\bf\thetitle}\par\medskip}   
\vglue 0.05truein 

{\parskip=0pt\leftskip 0pt plus 1fil\def\\{\par}{\sc\theauthors}
\par\medskip}%
 
\vglue 0.03truein 

{\small\leftskip 25pt\rightskip 25pt{\bf Abstract}\stdspace\theabstract

{\bf AMS Classification}\stdspace\theprimaryclass
\ifx\thesecondaryclass\relax\else; \thesecondaryclass\fi\par
{\bf Keywords}\stdspace \thekeywords\par}\vglue 7pt

}   

\font\phead=cmsl9 scaled 950
\font=cmsl9 scaled 1050
\font\pnum=cmbx10 scaled 913
\font\lnum=cmbx10 
\font\pfoot=cmsl9 scaled 950
\font=cmsl9 scaled 1050
\ifplaintex
\headline{\vbox to 0pt{\vskip -4.5mm\line{\small\phead\ifnum
\count0=\startpage ISSN 1464-8997 (on line)
1464-8989 (printed) \hfill {\pnum\folio}\else\ifodd\count0\def\\{ }%
\ifx\theshorttitle\relax\thetitle\else\theshorttitle\fi\hfill{\pnum\folio}
\else\def\\{ and }{\pnum\folio}\hfill\ifx\theshortauthors\relax\theauthors
\else\theshortauthors\fi\fi\fi}\vss}}
\footline{\vbox to 0pt{\vglue 0mm\line{\small\pfoot\ifnum\count0=\startpage
Published \publishdate:\qua\copyright\ \gtp\hfill\else
\gtm, Volume \thevolumenumber\ (\thevolumeyear)\hfill\fi}\vss
}}
\else
\makeatletter
\def\@oddhead{{\small\lhead\ifnum\count0=\startpage ISSN 1464-8997 (on line)
1464-8989 (printed) \hfill {\lnum\number\count0}\else\ifodd\count0
\def\\{ }\ifx\theshorttitle\relax \thetitle \else\theshorttitle\fi\hfill
{\lnum\number\count0}\else\def\\{ and }{\lnum\number\count0}
\hfill\ifx\theshortauthors\relax 
\theauthors\else\theshortauthors\fi\fi\fi}}\def\@evenhead{@oddhead}
\def\@oddfoot{\small\lfoot\ifnum\count0=\startpage Published \publishdate:\qua\copyright\ \gtp\hfill\else
\gtm, Volume \thevolumenumber\ (\thevolumeyear)\hfill\fi}
\def\@evenfoot{@oddfoot}
\makeatother
\fi

\let\maketitlepage\maketitlep
\let\makeshorttitle\maketitlepage
\let\maketitle\maketitlepage

\newwrite\gtoutfile
\long\gdef\makeheadfile{  
{\def\\{, }\def\s{ }
\immediate\openout\gtoutfile head.xxx
\immediate\write\gtoutfile{To: math@arxiv.org}
\immediate\write\gtoutfile{Subject: put OR rep NNNNN:ppppp}
\immediate\write\gtoutfile{--text follows this line--}
\immediate\write\gtoutfile{Proxy-for: \ifx\theasciiauthors\relax
\theauthors\else\theasciiauthors\fi\s<\ifx\theasciiemail\relax\theemail\else\theasciiemail\fi>}
\immediate\write\gtoutfile{\noexpand\\}
\immediate\write\gtoutfile{Authors: \ifx\theasciiauthors\relax
\theauthors\else\theasciiauthors\fi}
{\def\\{ }\immediate\write\gtoutfile{Title: \ifx\theasciititle\relax
\thetitle\else\theasciititle\fi}}
\immediate\write\gtoutfile{Subj-class: GT or SG, GR etc}
\immediate\write\gtoutfile{MSC-class: \theprimaryclass\ifx\thesecondaryclass\relax\else, \thesecondaryclass\fi}
\immediate\write\gtoutfile{Journal-ref: Geom. Topol. Monogr. \thevolumenumber\s
(\thevolumeyear) \startpage-\finishpage}
\immediate\write\gtoutfile{Comments: Published by Geometry and Topology Monographs at}
\immediate\write\gtoutfile{\s\s\s  http://www.maths.warwick.ac.uk/gt/GTMon\thevolumenumber/paper\thepapernumber.abs.html}
\immediate\write\gtoutfile{\noexpand\\}
\immediate\write\gtoutfile{}
\ifx\theasciiabstract\relax
\immediate\write\gtoutfile{\theabstract}\else
\immediate\write\gtoutfile{\theasciiabstract}\fi
\immediate\write\gtoutfile{}
\immediate\write\gtoutfile{\noexpand\\}
\immediate\write\gtoutfile{}
\immediate\closeout\gtoutfile}}  

\def\maketitlepage{\maketitlep\makeheadfile}
\let\makeshorttitle\maketitlepage
\let\maketitle\maketitlepage

\volumenumber{4}
\volumename{Invariants of knots and 3-manifolds (Kyoto 2001)}
\volumeyear{2002}
\papernumber{21}
\pagenumbers{313}{335}
\received{8 November 2002}
\revised{17 October 2003}
\accepted{1 November 2003}
\published{13 November 2003}

\usepackage{graphicx,amsmath}

\numberwithin{figure}{section}

\newtheorem{theorem}{Theorem}[section]
\newtheorem{lemma}[theorem]{Lemma}
\newtheorem{corollary}[theorem]{Corollary}
\newtheorem{proposition}[theorem]{Proposition}

\newtheorem{condition}[theorem]{Condition}

\newtheorem{conjecture}[theorem]{Conjecture}
\theoremstyle{definition}

\newtheorem{definition}[theorem]{Definition}
\newtheorem{problem}[theorem]{Problem}
\newtheorem{example}[theorem]{Example}

\begin{document}
\title{Skein module deformations of\\elementary moves on links}
\author{J\'ozef H Przytycki}
\asciiauthors{Jozef H Przytycki}
\coverauthors{J\noexpand\'ozef H Przytycki}
\address{Department of Mathematics,
George Washington University\\Washington, DC 20052,  USA}
\email{przytyck@gwu.edu}

\begin{abstract}
This paper is based on my talks (``Skein modules with a cubic skein
 relation: properties and speculations'' and ``Symplectic structure on
 colorings, Lagrangian tangles and its applications") given in Kyoto
 (RIMS), September 11 and September 18 respectively,
 2001.  The first three sections closely follow the talks: starting
 from elementary moves on links and ending on applications to
 unknotting number motivated by a skein module deformation of a
 3-move.  The theory of skein modules is outlined in the
 problem section of these proceedings (see \cite{H-P,P-3,Tu}).

In the first section we make the point that despite its long history,
knot theory has many elementary problems that are still open.  We
discuss several of them starting from the Nakanishi's 4-move
conjecture.  In the second section we introduce the idea of {\it
Lagrangian tangles} and we show how to apply them to elementary moves
and to rotors.  In the third section we apply $(2,2)$-moves and a
skein module deformation of a 3-move to approximate unknotting numbers
of knots.  In the fourth section we introduce the Burnside groups of
links and use these invariants to resolve several problems stated in
Section 1.
\end{abstract}

\asciiabstract{This paper is based on my talks (`Skein modules with a 
cubic skein relation: properties and speculations' and `Symplectic
 structure on colorings, Lagrangian tangles and its applications')
 given in Kyoto (RIMS), September 11 and September 18 respectively,
 2001.  The first three sections closely follow the talks: starting
 from elementary moves on links and ending on applications to
 unknotting number motivated by a skein module deformation of a
 3-move.  The theory of skein modules is outlined in the problem
 section of these proceedings.

In the first section we make the point that despite its long history,
knot theory has many elementary problems that are still open.  We
discuss several of them starting from the Nakanishi's 4-move
conjecture.  In the second section we introduce the idea of
Lagrangian tangles and we show how to apply them to elementary moves
and to rotors.  In the third section we apply (2,2)-moves and a
skein module deformation of a 3-move to approximate unknotting numbers
of knots.  In the fourth section we introduce the Burnside groups of
links and use these invariants to resolve several problems stated in
section 1.}


\maketitle

\section{Elementary moves on links: history and conjectures}

Knot theory is more than two hundred years old; the first scientists
who considered knots as mathematical objects were A Vandermonde (1771) and
C\,F Gauss (1794). However, despite the impressive growth of the theory,
there are simply formulated yet fundamental questions, to which we do not
know answers.
These problems are not just interesting puzzles but they lead  
to very interesting theory (structure).

The oldest of such problems is the Nakanishi 4-move conjecture,
formulated in 1979 \cite{Kir,Nak-1,P-1}.
Recall that an $n$-move on a link is a local change of the 
link illustrated in Figure 1.1.
In our convention the part of the link outside of the disk,
in which the move takes part, is unchanged.

Two unoriented links are said to be $n$-move equivalent 
if there is a sequence of
$\pm n$-moves which converts one link to the other.

\begin{figure}[ht!]
\centerline{\includegraphics[height=1.4cm]{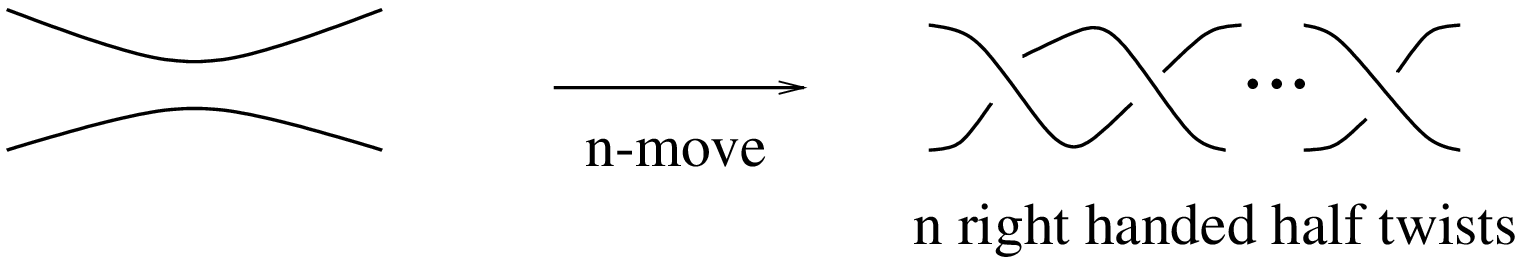}}
\caption{$n$-move}
\end{figure}

\begin{conjecture}[Nakanishi, 1979]\label{1.1}
Every knot is $4$-move equivalent to the trivial knot.
\end{conjecture}

The conjecture holds for closed 3-braids and 2-algebraic links (i.e.\
algebraic links in the Conway sense) \cite{Kir}. 
Because the concept of algebraic links \cite{Co,B-S}
 and its generalizations will
be used often later in the paper, 
let us recall the definition \cite{P-Ts}.

\begin{definition}\label{1.2}$\phantom{99}$
\begin{enumerate}
\item[(i)] The set of $n$-algebraic tangles is the smallest family
of $n$-tangles which satisfies:

(0) Any $n$-tangle with 0 or 1 crossing is $n$-algebraic.

(1) If $A$ and $B$ are $n$-algebraic tangles then $r^i(A)*r^j(B)$
is $n$-algebraic; $r$ denotes here the rotation of a tangle 
along the $z$-axis by the angle $\frac{2\pi}{2n}$, 
and * denotes (horizontal) composition of tangles.

\item[(ii)] If in the condition (1), $B$ is restricted
to tangles with no more than $k$ crossings, we obtain
the family of $(n,k)$-algebraic tangles.

\item[(iii)]
If a link, $L$, is obtained from an $(n,k)$-algebraic
tangle (resp. $n$-algebraic tangle)
by connecting its endpoints 
without introducing any new crossings then $L$ is called
an $(n,k)$-algebraic (resp. $n$-algebraic) link. 
\end{enumerate}
\end{definition} 
Two examples of 2 and 4-algebraic tangles are shown in Figure 1.2.

\begin{figure}[ht!]
\centerline{\includegraphics[height=1.6cm]{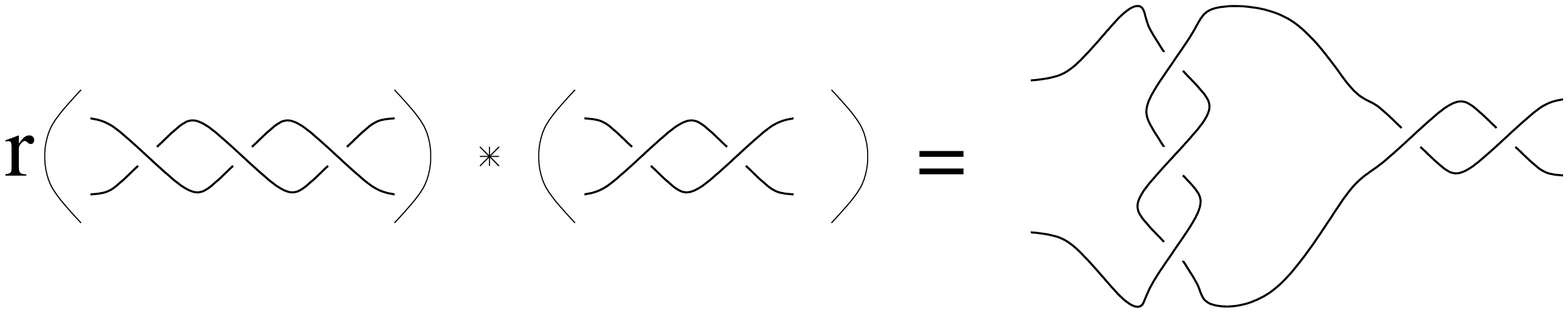}}
\vspace{3mm}

\centerline{\includegraphics[height=2.1cm]{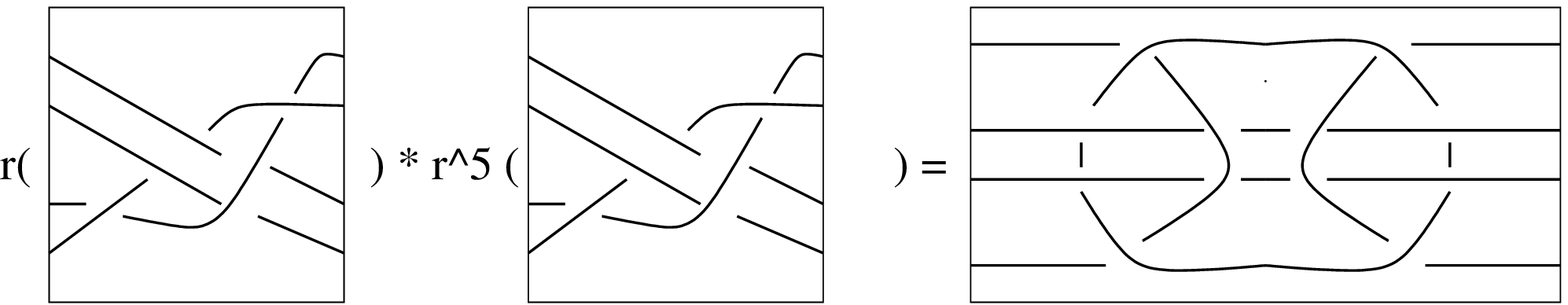}}
\nocolon\caption{}
\end{figure}

The inductive character of our definition allows us to show:
\begin{proposition}{\rm\cite{P-5}}\label{1.3}
\begin{enumerate}
\item[\rm(a)]
Every 2-algebraic tangle without a closed component can
be reduced by 4-moves to one of the following six 2-tangles,
Figure  1.3.
\item[\rm(b)]
Every 2-algebraic knot can be reduced by 4-moves to the trivial knot.
\end{enumerate}
\end{proposition}

\begin{figure}[ht!]
\centerline{\includegraphics[height=1.6cm]{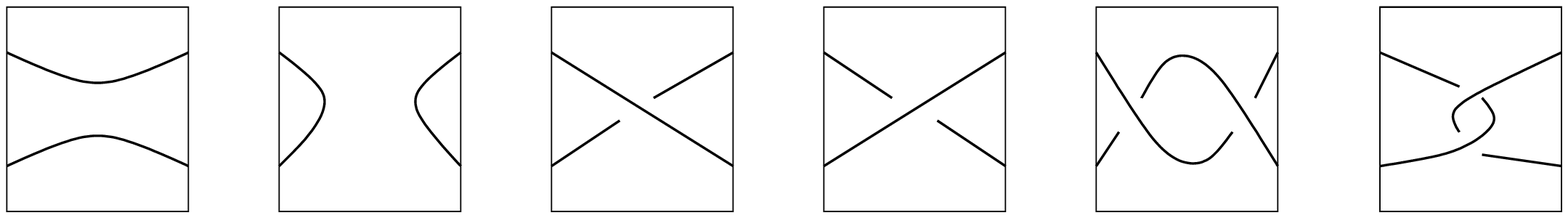}}
\nocolon\caption{}
\end{figure}

In 1994 Nakanishi began to suspect that the $(2,1)$-cable
of the trefoil knot cannot be simplified by 4-moves \cite{Kir}. 
However N.Askitas was able to simplify this knot \cite{Ask}. 
Askitas, in turn,
suspects the $(2,1)$-cable of the figure eight knot as the
simplest counter-example to the Nakanishi 4-move conjecture. 

Not every link can be reduced to a trivial link by 4-moves,
in particular the linking matrix modulo 2 is preserved by 4-moves.
Furthermore Nakanishi demonstrated that the 
Borromean rings cannot be reduced to the trivial link 
of 3-components \cite{Nak-2}.

Kawauchi expressed the question for links as follows:
\begin{problem}\cite{Kir}
\begin{enumerate}
\item[(i)]
Is it true
that if two links are link-homotopic 
then they are 4-move equivalent?
\item[(ii)] In particular, is it true
that every 2-component link is 4-move equivalent to the trivial link
of two components or to the Hopf link?
\end{enumerate}
\end{problem}
We can use an inductive argument to show that
any 2-component 2-algebraic link is 4-move equivalent 
to a trivial link or a Hopf link.
We also have proved that the answer to Kawauchi's question is affirmative 
for closed 3-braids \cite{P-5}. 

Nakanishi identified the ``half" 
2-cabling of the Whitehead link, ${\cal W}$,
 Figure 1.5,
as the simplest link which he could not reduce by 4-moves but which
is link homotopy equivalent to the trivial link\footnote{In June of 2002
we have shown that this example cannot, in fact, be reduced by 4-moves,
 \cite{D-P-2}. We discuss this result in Section 4.}.
\addtocounter{figure}{1}

\begin{figure}[ht!]
\centerline{\includegraphics[height=3.5cm]{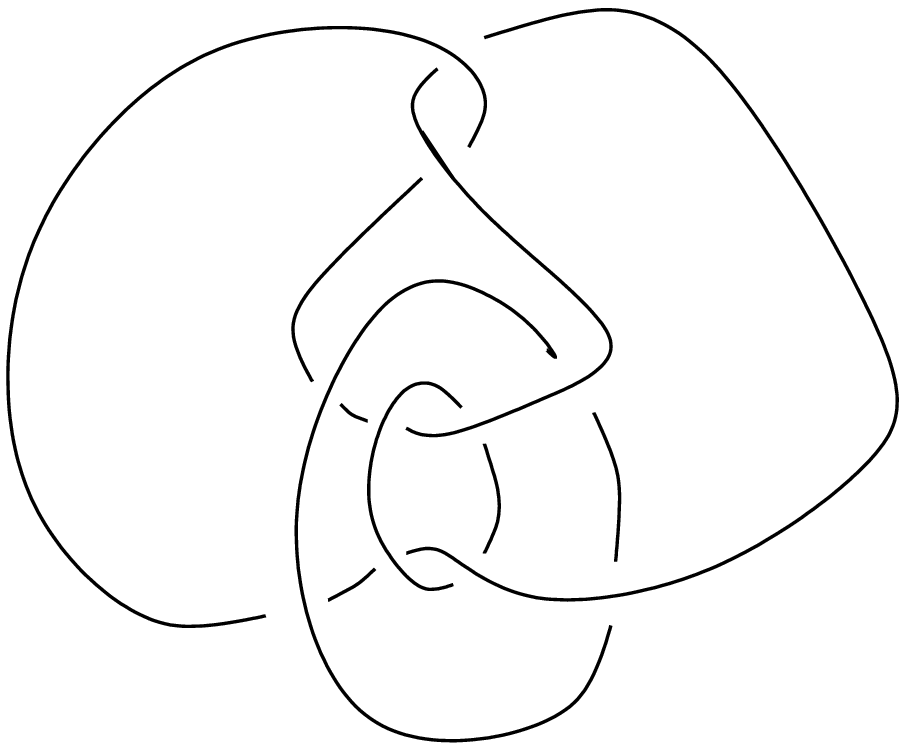}}
\nocolon\caption{}
\end{figure}

The second oldest and the most studied elementary moves question 
is the Montesinos-Nakanishi 3-move conjecture. 
\begin{conjecture}\label{1.5}
Any link is 3-move equivalent to a trivial link.
\end{conjecture}
Nakanishi first considered the conjecture in 1981. J.~Montesinos
analyzed 3-moves before, in connection with 3-fold dihedral branch
coverings \cite{Mon}.

The conjecture has been proved to be valid for several classes of links:
by Y.Nakanishi for links up to 10 crossings and Montesinos links,
by J.Przytycki for links up to 11 crossings, 2-algebraic
links and closed 3-braids, by Q.Chen for links up to 12 crossings,
closed 4-braids and closed 5-braids with the exception of the class
of $\hat\gamma$ -- the square of the center of the fifth braid group
($\gamma = (\sigma_1 \sigma_2 \sigma_3 \sigma_4)^{10}$),
by Przytycki and T.Tsukamoto
for 3-algebraic links (including 3-bridge links), and by Tsukamoto
for $(4,5)$-algebraic links (including 4-bridge links),
\cite{Kir,Che,P-Ts,Tsu}.

Nakanishi presented in 1994, an example which he
could not reduce by 3-moves: the 2-parallel of the Borromean rings,
$L_{2BR}$. The link $L_{2BR}$ has 24 crossings and $\hat\gamma$ can be
reduced to a link of 20 crossings (Figure 1.6) which is the smallest
link not reduced yet by 3-moves\footnote{We proved in \cite{D-P-1} 
that neither
$\hat\gamma$ nor $L_{2BR}$ can be reduced by 3-moves to a trivial
link; see Section 4.}.

\begin{figure}[ht!]
\centerline{\includegraphics[height=5.1cm]{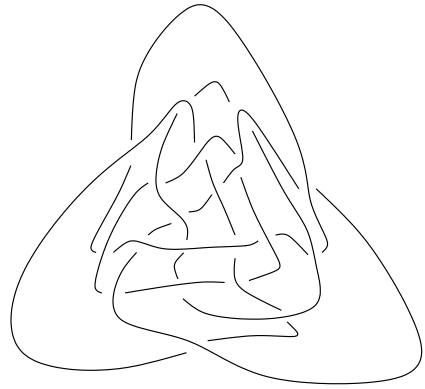}}
\caption{$L_{2BR}$ and the Chen's link}
\end{figure}

Not every link can be simplified using
$5$-moves, but every $5$-move is a combination of more
delicate  $(2,2)$-moves, or, more precisely of a $(2,2)$-move
(\lower .5cm \hbox{\includegraphics[height=1cm]{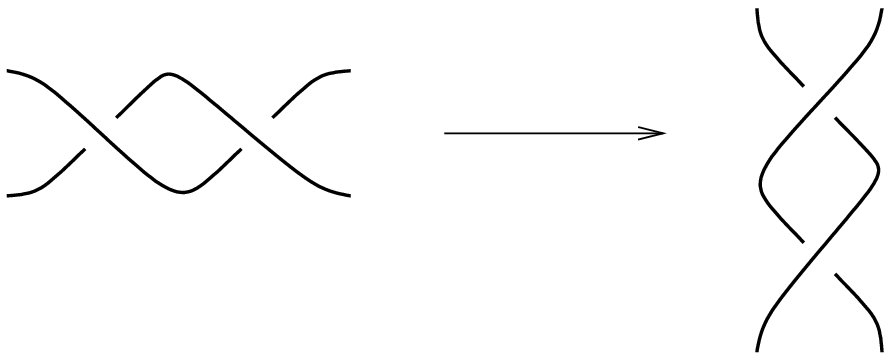}}) and
its mirror image $(-2,-2)$-move 
(\lower .4cm \hbox{\includegraphics[height=1cm]{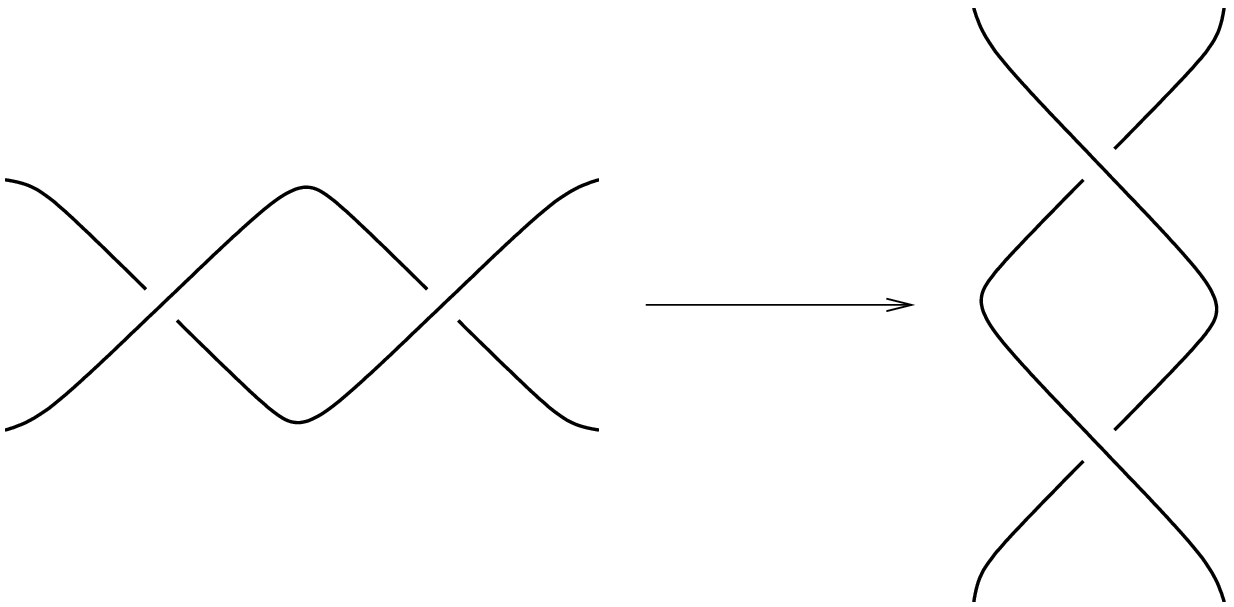}}).

\begin{figure}[ht!]
\centerline{\includegraphics[height=5cm]{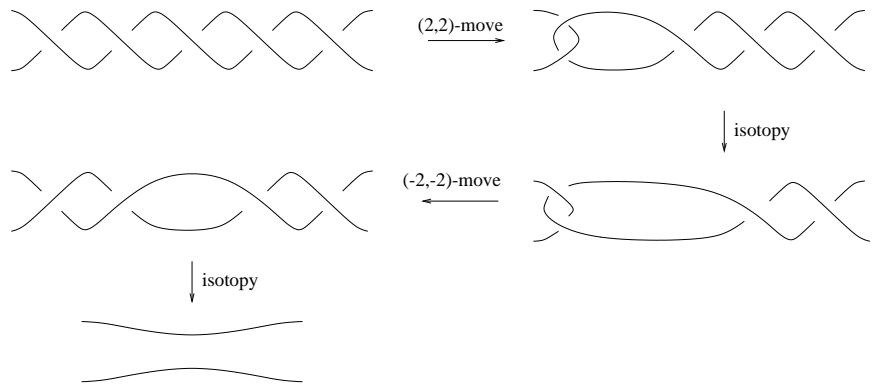}}
\caption{$5$-move as a combination of $\pm(2,2)$-moves}
\end{figure}

\begin{conjecture}[Harikae, Nakanishi, Uchida, 1992]\label{1.6}

Every link is $(2,2)$-move equivalent to a trivial link.
\end{conjecture}
Conjecture 1.6 has been established for several classes of links.
\begin{lemma}{\rm \cite{P-5}}\qua\label{1.7}
Every 2-algebraic link is $(2,2)$-move equivalent to a trivial link.
\end{lemma}
\begin{lemma}{\rm\cite{P-5}}\qua\label{1.8}
Every link, up to 9 crossings, is $(2,2)$-move equivalent either
to a trivial link or to $9_{40}$ or $9_{49}$ knots\footnote{
We showed
in \cite{D-P-2} that the knots $9_{40}$ and $9_{49}$ are
 not $(2,2)$-move equivalent to trivial links; see Section 4.}
, or 
to their mirror images $\bar 9_{40}$ and $\bar 9_{49}$.
\end{lemma}
\begin{lemma}{\rm\cite{P-5}}\qua\label{1.9}
Every closed 3-braid is $(2,2)$-move equivalent to a trivial link
or to the closure of the braids $(\sigma_1^2 \sigma_2^{-1})^3$ or 
$\sigma_1^2 \sigma_2^2\sigma_1^{-2}\sigma_2^2\sigma_1^2\sigma_2^{-2}$,
or their mirror images.
\end{lemma}
Notice that  the knot $9_{49}$ and the closure of $(\sigma_1^2\sigma_2^{-1})^3$
(i.e.\ the link $9^2_{40}$ in Rolfsen's notation \cite{Rol}) are
related by a $(2,2)$-move; Figure 1.8. Similarly, the knot $9_{40}$
and the closure of the 3-braid 
$\sigma_1^2 \sigma_2^2\sigma_1^{-2}\sigma_2^2\sigma_1^2\sigma_2^{-2}$
are $(2,2)$-move equivalent \cite{P-5}.

\begin{figure}[ht!]
\centerline{\includegraphics[height=4.3cm]{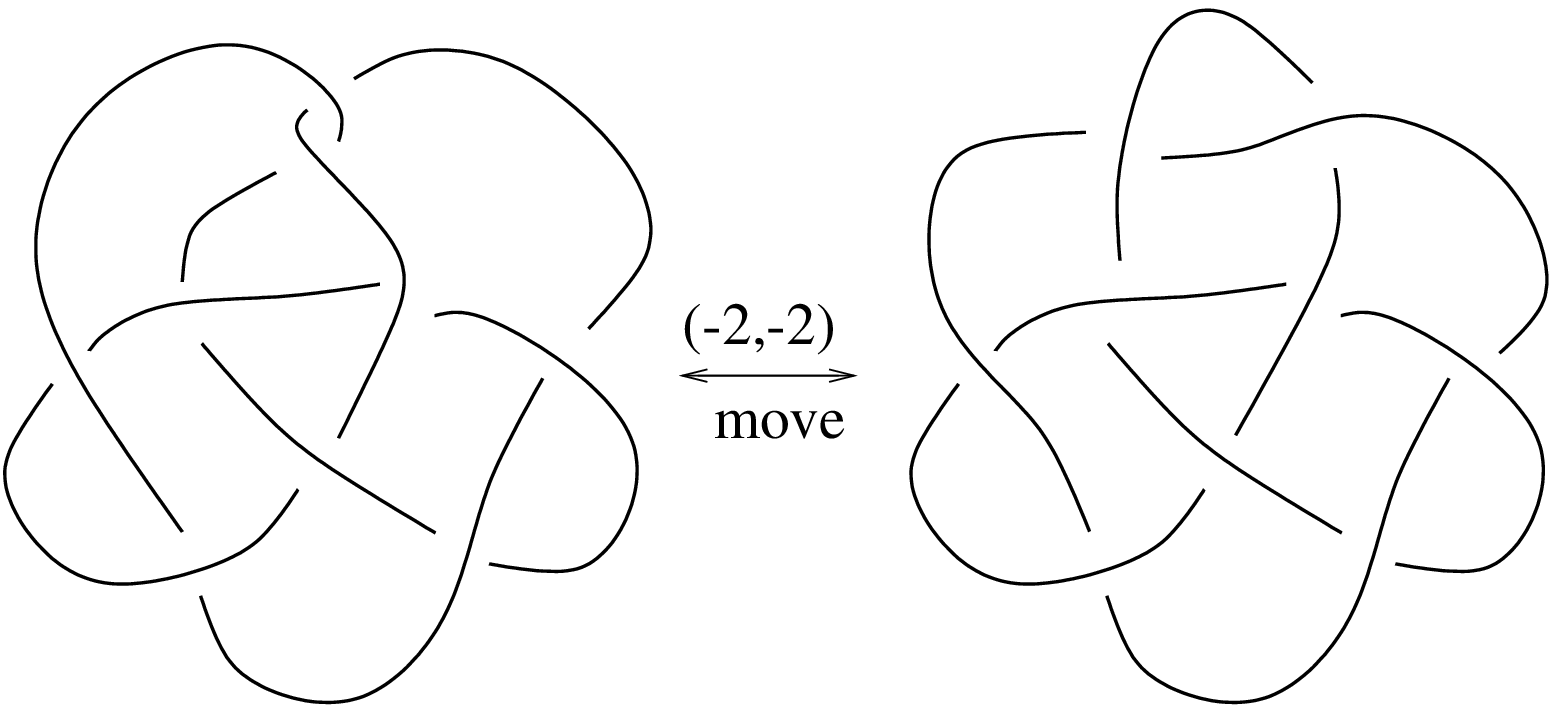}}
\nocolon\caption{}
\end{figure}

Similarly as in the case of $5$-moves, not
every link can be reduced via
$7$-moves to a trivial link. The $7$-move is however a combination of
$(2,3)$-moves\footnote{To be precise, a $7$-move is a combination
of a $(-3,-2)$- and a $(2,3)$-move; compare Figures 1.9 and 1.7.}
which might be sufficient for a reduction. We say that two links are
$(2,3)$-move equivalent if there is a sequence of
$\pm(2,3)$-moves and their inverses, $\pm(3,2)$-moves,
 which converts one to the other.
\begin{problem}(Kir)\qua\label{1.10}
Is every link $(2,3)$-move equivalent to a trivial link?
\end{problem}

\begin{figure}[ht!]
\centerline{\includegraphics[height=3cm]{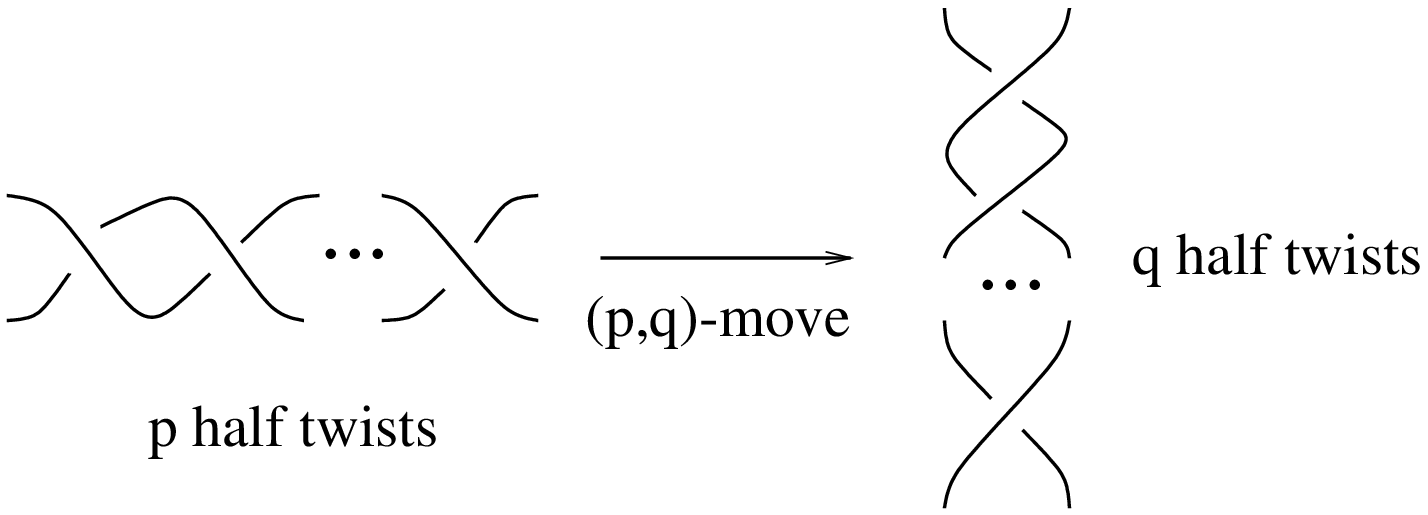}}
\nocolon\caption{}
\end{figure}

The answer is affirmative for 2-algebraic links.
Furthermore every 2-algebraic tangle can be reduced to
one of the 8 ``basic" algebraic tangles (with possible
additional trivial components) presented in Figure 1.10.

\begin{figure}[ht!]
\centerline{\includegraphics[height=3.6cm]{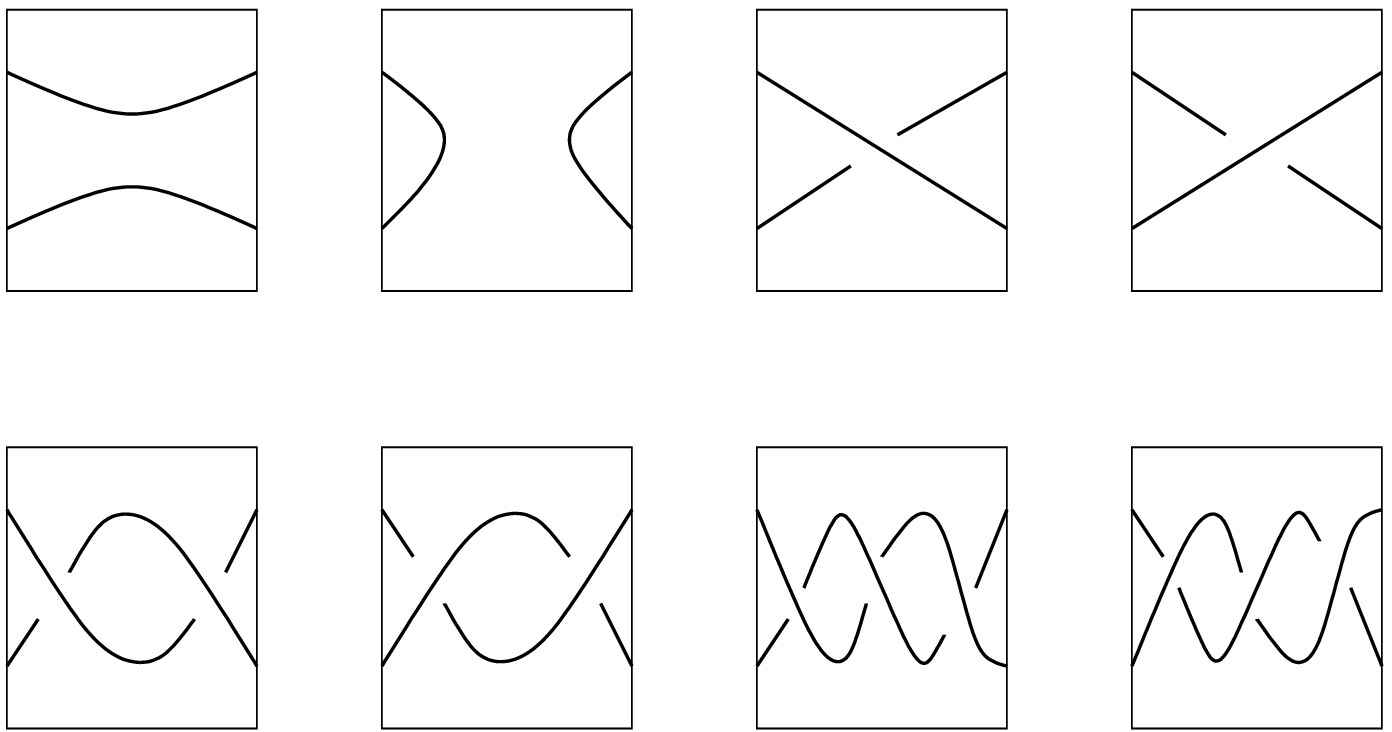}}
\nocolon\caption{}
\end{figure}

Six months ago, in March of 2000, presenting the talk
``Open problems in knot theory that everyone can try to solve"
at the conference in Cuautitlan
I tried to extend the range
of classical unknotting moves, so I proposed the following question.
\begin{problem}\label{1.11}
Can every link be reduced to a trivial link by $(2,5)$ and $(4,-3)$ moves, 
their inverses and mirror images? 
\end{problem}

As before the question was motivated by the observation that the answer
is positive for 2-algebraic links. In the proof we noticed that
every 2-algebraic tangle can be reduced to
one of the 12 ``basic" algebraic tangles (with possible
additional trivial components) presented in Figure 1.11.

\begin{figure}[ht!]
\centerline{\includegraphics[height=4.6cm]{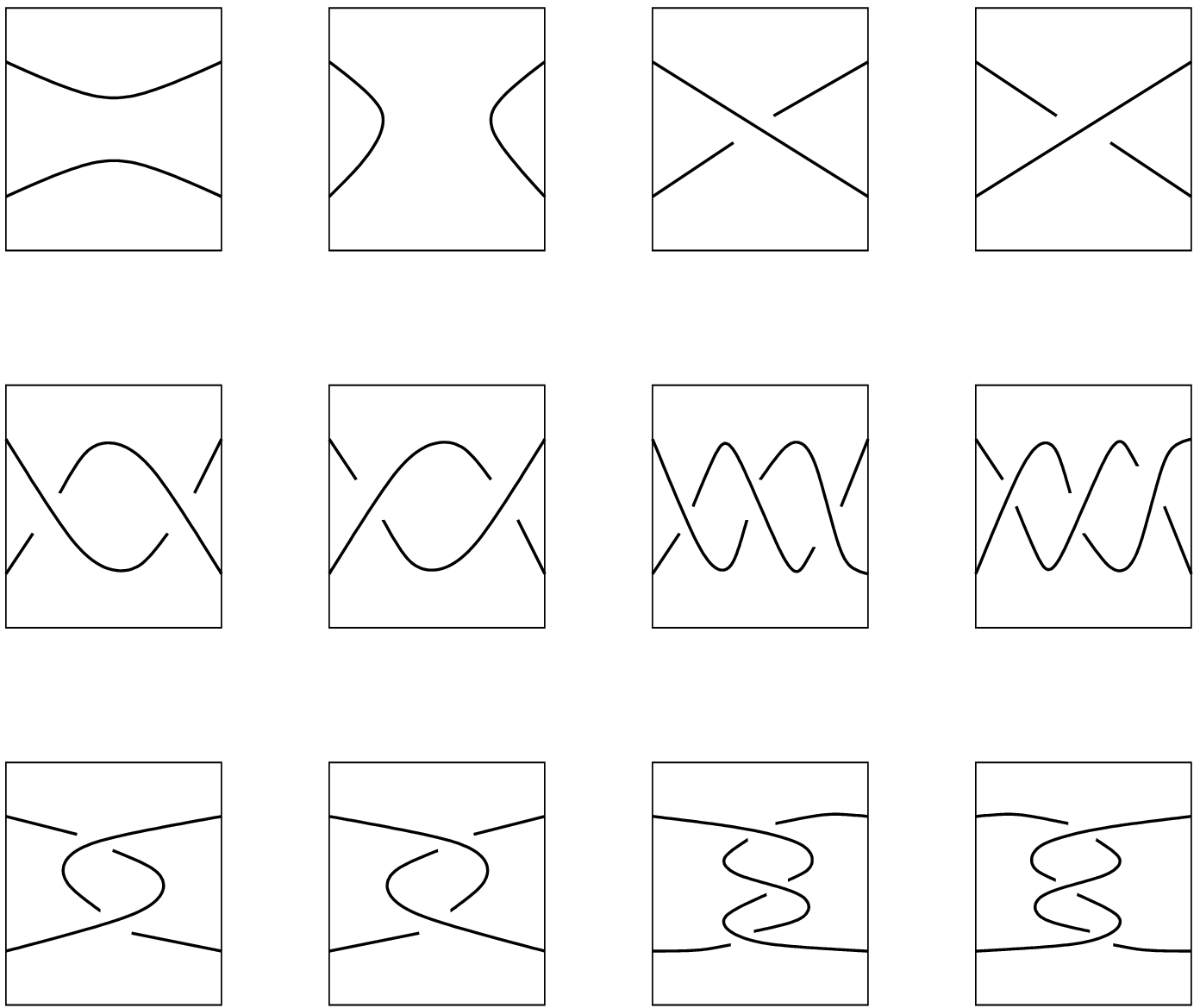}}
\nocolon\caption{}
\end{figure}

We can ask: how to generalize further our problem? 
Notice that $(2,5)$ and $(4,-3)$ moves
are equivalent to $2+ \frac{1}{5} =\frac{11}{5}$ and 
$4+ \frac{1}{-3} =\frac{11}{3}$ rational moves, respectively.
Figure 1.12 illustrates the fact
that an $(s,q)$-move is equivalent to the rational $s+\frac{1}{q}=
\frac{sq+1}{q}$-move.

\begin{figure}[ht!]
\centerline{\includegraphics[width=.9\hsize]{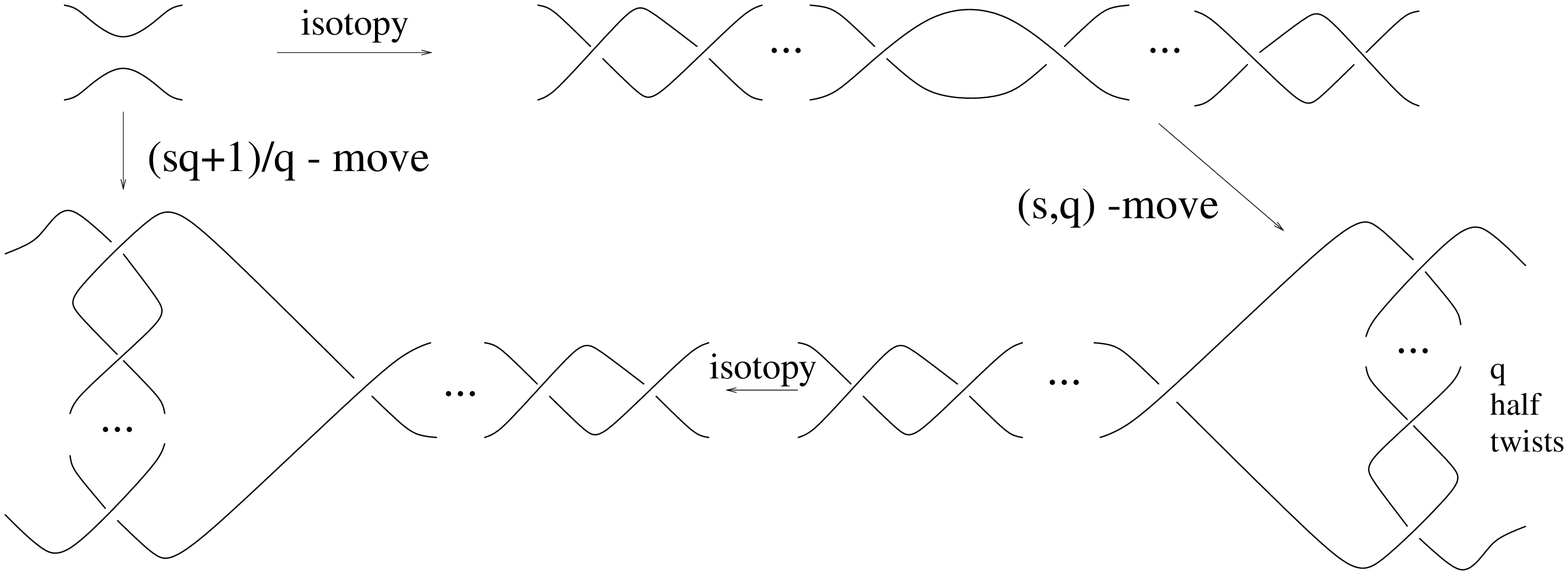}}
\caption{$(s,q)$-move as $\frac{sq+1}{q}$-move}
\end{figure}

These suggest the following extension of the problem using 
rational moves\footnote{In each $\frac{p}{q}$-rational move, 
a $0$-tangle is replaced by a rational $\frac{p}{q}$-tangle.}.
\begin{problem}\label{1.12}
Let $p$ be a fixed prime number, then
\begin{enumerate}
\item[(i)] Is it true that every link can
be reduced to a trivial link
by rational $\frac{p}{q}$-moves ($q$ any integer).\footnote{We proved 
in \cite{D-P-2} that the answer is negative: for $p\geq 5$ the
closure of the 3-braid $(\sigma_1\sigma_2)^6$ cannot be reduced; compare
Section 4.}
\item[(ii)] Is there a function
$f(n,p)$ such that any $n$-tangle can be reduced to one of
``basic" $f(n,p)$ $n$-tangles (allowing additional trivial components)
by rational $\frac{p}{q}$-moves.
\end{enumerate}
\end{problem}
The method of ``Lagrangian tangles"\cite{D-J-P} has allowed us to prove 
that $f(n,p) \geq \Pi_{i=1}^{n-1}(p^i+1)$. We will discuss the method in
the second section.

\section{Symplectic structures, Lagrangian tangles and Rotors}\label{2} 

I will illustrate here how the symplectic structure can be used to answer
the ``classical" question about homology of double branched covers 
(or Fox colorings). 
I will describe how symplectic structures came into the picture and
finally use them to analyze how homology of double branched covers 
changes after rotation \cite{A-P-R,D-J-P,P-4}.
\begin{definition}[Fox colorings]\label{2.1}$\phantom{99}$
\begin{enumerate}
\item [(i)]
We say that a link (or a tangle) diagram is k-colored if every
arc is colored by one of the numbers $0,1,...,k-1$ (forming a
group $Z_k$) in such a way that
at each crossing the sum of the colors of the undercrossings is equal
to twice the color of the overcrossing modulo $k$; see Figure 2.1.
\item [(ii)]
The set of $k$-colorings forms an abelian group, denoted by $Col_k(D)$.
The cardinality of the group will be denoted by $col_k(D)$.
For an $n$-tangle $T$ each Fox $k$-coloring of $T$ yields a
coloring of boundary points of $T$ and we have the homomorphism
$\psi :Col_k(T) \rightarrow Z_k^{2n}$.
\end{enumerate}
\end{definition}.

\begin{figure}[ht!]
\centerline{\includegraphics[height=1.6cm]{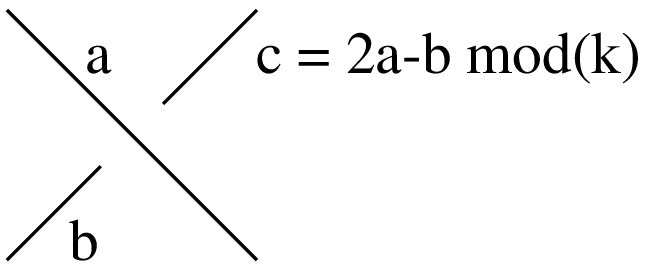}}
\nocolon\caption{}
\end{figure}

It is a pleasant exercise to show that $Col_k(D)$ is unchanged
by Reidemeister moves and by $k$-moves (Figure 2.2).

\begin{figure}[ht!]
\centerline{\includegraphics[height=1.6cm]{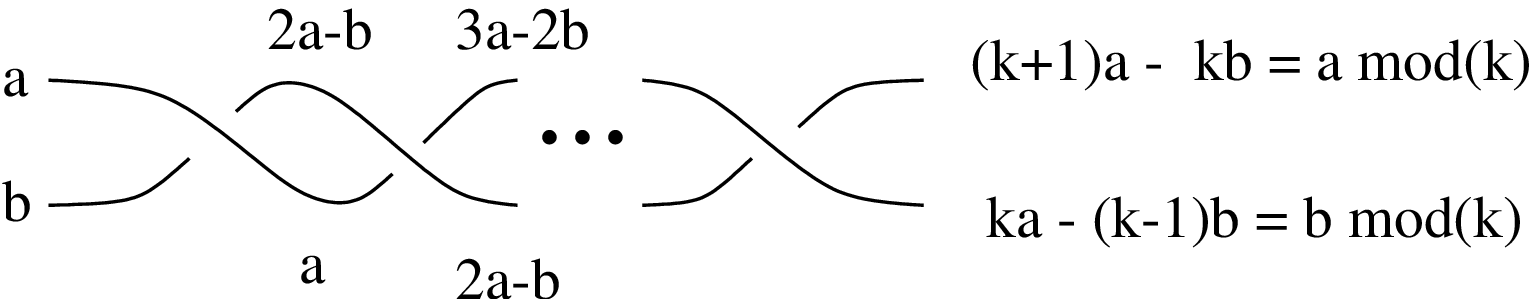}}
\nocolon\caption{}
\end{figure}

It is known that the module $Col_k(T)= H^1(M^{(2)}_T,Z_k)\oplus Z_k$ where
$M^{(2)}_T$ denotes the double branched cover of $D^3$ with the
tangle, $T$, as the branch set (compare \cite{P-2}).

I was curious for a long time how to characterize images 
$\psi(Col_k(T))$. I was sure that such a description could be used 
to analyze elementary moves on links\footnote{
At the time I was giving my talks in Kyoto I had high hopes
that the Fox 3-colorings
are the only obstructions to 3-move equivalence of links. I extended this
hope to tangles. From May of 2000 I was very much involved in an analysis of
the structure of the space of 3-colorings of tangles, $Col_3(T)$.
In particular, we had
two invariants of 3-move equivalence classes of tangles: the space $Col_3(T)$
and its image under the map $\psi: Col_3(T) \to Z_3^{2n}$, where
$Z_3^{2n}$ represents the space of colorings of $2n$ boundary points of $T$.
By interpreting tangles as Lagrangian subspaces in the symplectic space
$Z_3^{2n-2}$ we were able to prove that there are exactly $\Pi_{i=1}^{n-1}
(3^i+1)$ different images $\psi(Col_3(T))$. Obstructions to the
Montesinos-Nakanishi conjecture, constructed in February of 2002 (see
Section 4),
show however that there are more 3-move equivalence
classes of $n$-tangles than predicted by 3-coloring invariants \cite{D-P-1}.
To understand these classes is the major task.}. 
Below I describe the solution in
which spaces $\psi(Col_k(T))$ are identified as Lagrangians in
some symplectic spaces \cite{D-J-P}. 

Consider $2n$ points on a circle (or a square) and a
field ${\bf Z}_p$ of $p$-colorings of a point. The colorings of
$2n$ points form ${\bf Z}_{p}^{2n}$ linear space. Let $e_1,\ldots
, e_{2n}$ be its basis,

\begin{figure}[ht!]
\centerline{\includegraphics[height=1.3cm]{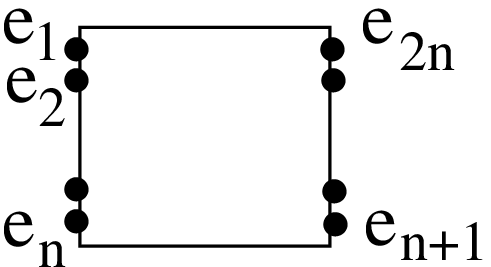}}
\nocolon\caption{}
\end{figure}

$e_i=(0,\ldots , 1,\ldots ,0)$, where 1
occurs in the $i$-th position. Let ${\bf Z}_{p}^{2n-1}\subset{\bf
Z}_p^{2n}$ be the subspace of vectors $\sum a_i e_i$ satisfying $\sum (-1)^i
a_i=0$ (alternating condition).
Consider the basis $f_1,\ldots , f_{2n-1}$ of ${\bf
Z}_{p}^{2n-1}$ where $f_k=e_k+e_{k+1}$. Consider a skew-symmetric
form $\phi$ on ${\bf Z}_{p}^{2n-1}$ of nullity $1$ given by the
matrix

$$\phi = \left( \begin{array}{cccc} 
0 & 1 & 0 &\ldots \\
 -1 & 0 & 1 &\ldots \\ 
\ldots & \ldots & \ldots & \ldots \\ 
 0 & \ldots & -1 & 0
\end{array} \right)$$

that is
 \vspace{1mm}
   \renewcommand{\arraystretch}{1.5}
   $$\phi (f_i, f_j) =\left\{
   \begin{array}{lr}
   0 &
   {\rm if}\ |j-i|\neq 1 \\
  1 &  {\rm if\ }\ j=i+1\\ -1 &
   {\rm if}\ j=i-1. \\
   \end{array}
  \right .$$
   \par
   \vspace{2mm}

Notice that the vector $e_1+ e_2+ \ldots + e_{2n}$
($=f_1+f_3+\ldots + f_{2n-1}=f_2+f_4+\ldots + f_{2n})$
is $\phi$-orthogonal to any other vector.
If we consider ${\bf Z}_{p}^{2n-2}={\bf Z}_{p}^{2n-1}/{\bf
Z}_{p}$, where the subspace ${\bf Z}_{p}$ is generated by
$e_1+\ldots + e_{2n}$, that is, ${\bf Z}_{p}$ consists
of monochromatic (i.e.\ trivial) colorings,
then $\phi$ descends to the symplectic form $\hat\phi$ on ${\bf
Z}_{p}^{2n-2}$. Now we can analyze isotropic subspaces of
$({\bf Z}_{p}^{2n-2},\hat\phi)$, that is subspaces on
which  $\hat\phi$ is $0$
($W\subset {\bf Z}_{p}^{2n-2}, \phi (w_1,w_2)=0$ for $w_1,w_2\in W$).
The maximal isotropic ($(n-1)$-dimensional)
subspaces of ${\bf Z}_{p}^{2n-2}$ are
called Lagrangian subspaces (or maximal totally degenerated subspaces)
and there are $\prod_{i=1}^{n-1}(p^i+1)$ of them.

We have $\psi :Col_p T\rightarrow{\bf Z}_{p}^{2n}$. Our local
condition on Fox colorings (Figure 2.1) guarantees that for any tangle $T$,
$\psi (Col_p T)\subset {\bf Z}_{p}^{2n-1}$. Furthermore,
the space of trivial colorings, ${\bf Z}_{p}$, always lays in $Col_p T$.
Thus $\psi$ descents to $\hat\psi :Col_p T/{\bf Z}_{p} \rightarrow
{\bf Z}_{p}^{2n-2}={\bf Z}_{p}^{2n-1}/{\bf Z}_{p}$. Now we have
the fundamental question:
Which subspaces of
${\bf Z}_{p}^{2n-2}$ are yielded by $n$-tangles?
We answer this question below.
\begin{theorem}{\rm\cite{D-J-P}}\qua\label{2.2}
$\hat\psi (Col_p T/{\bf Z}_{p})$ is a Lagrangian
subspace of ${\bf Z}_{p}^{2n-2}$ with the symplectic form $\hat\phi$.
\end{theorem}
The natural question would be whether every Lagrangian subspace
can be realized by a tangle. The answer is negative for
$p=2$ and positive for $p>2$ \cite{D-J-P}.
As a corollary we obtain a fact which was considered
to be difficult before, even for 2-tangles.
\begin{corollary}\label{2.3}
For any $p$-coloring of a tangle boundary satisfying
the alternating property (i.e.\ an element
of ${\bf Z}_{p}^{2n-1}$) there is an $n$-tangle and
its $p$-coloring yielding the given coloring on the boundary.
In other words: ${\bf Z}_{p}^{2n-1} = \bigcup_T \psi_T(Col_p(T))$.
Furthermore, the space $\psi_T(Col_p(T))$ is $n$-dimensional.
\end{corollary}


A few months ago we have found a rather spectacular application of the 
``Lagrangian tangles" method to rotors of links \cite{D-J-P}.
We have given criteria when $n$-rotation preserves the space of Fox
$p$-colorings, $Col_p(L)$.
\begin{definition}\cite{A-P-R}\qua\label{2.5}
Consider an $n$-tangle, that is a part of the link diagram, 
placed in the regular $n$-gon with 2n boundary points 
(n inputs and n outputs). We say that
this $n$-tangle is an $n$-rotor if it has a rotational symmetry, that is
the tangle is invariant with respect to rotation along $z$-axis by
the angle $\frac{2\pi}{n}$; see Figure 2.4.
\end{definition}

\begin{figure}[ht!]
\centerline{\includegraphics[height=4.7cm]{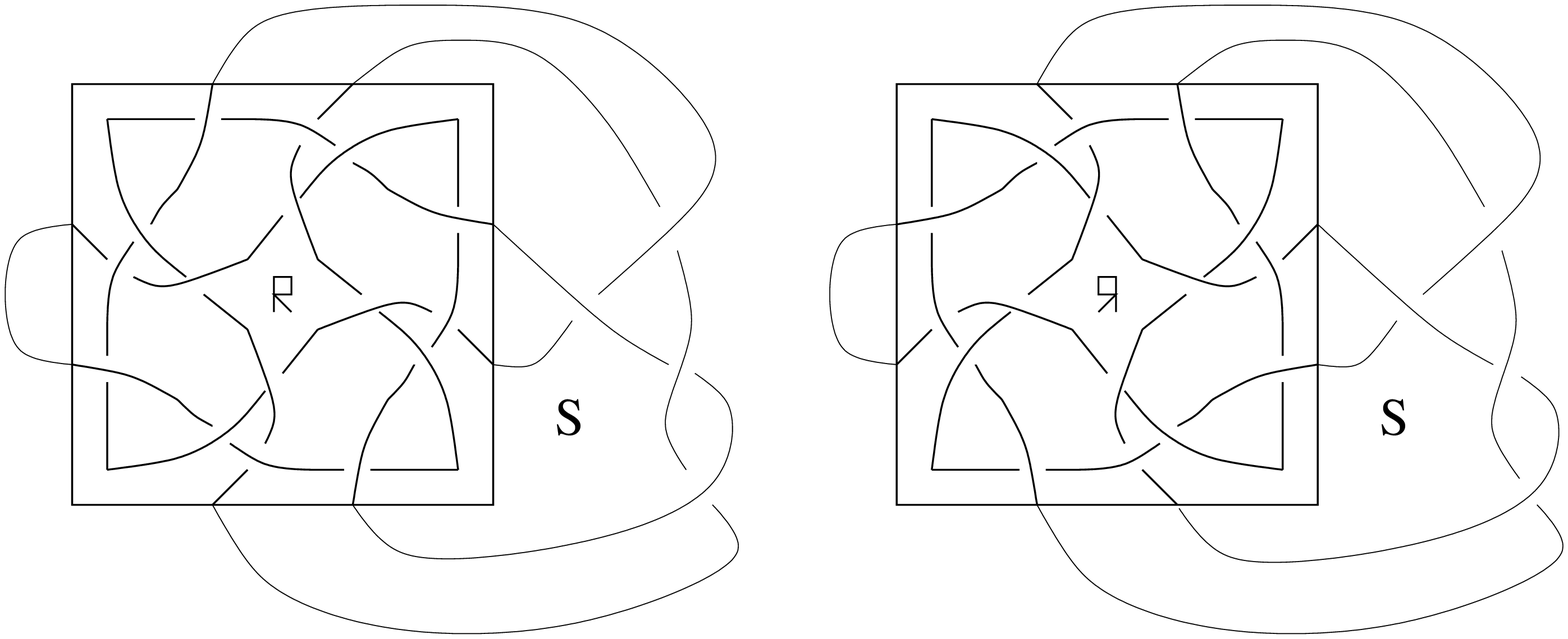}}
\nocolon\caption{}
\end{figure}

\begin{theorem}{\rm\cite{D-J-P}}\qua\label{2.6}
Let $L$ be a link diagram with an $n$-rotor part $R$. Let the rotant $\rho(L)$,
be obtained from $L$ by rotating $R$ around the $y$-axis by the angle $\pi$
and keeping the stator, $L-R$, unchanged. Assume that either $n=p$, 
where $p$ is a prime
number, or $n$ is co-prime to $p$ and 
such that there exists an $s$ with $p^s \equiv -1 \ mod\ n$.
Then the space $Col_p(L)$ is preserved by any $n$-rotation.
\end{theorem}
The technique we use in the proof of the theorem is to analyze
eigenspaces of the symplectic space $Z_p^{2n-2}$ with respect
to rotation. We obtain conditions under which the Lagrangian
subspace which is invariant under rotation is also invariant under
dihedral ``flype" $\rho$. Figure 2.4 illustrates a pair of 4-rotants
which have different space of Fox 5-colorings, $Z_5^2$ and $Z_5^3$,
respectively. One can also compute from Figure 2.4 that $ker\psi =Z_5$,
therefore there exists a nontrivial $5$-coloring which is equal to 0
on the boundary of $R$\footnote{It follows from the fact that 
any rotation preserves the determinant of $L$ which is equal to the 
order of $H_1(M^{(2)}_L,Z)$ (\cite{P-4}). Thus if a rotation changes 
the space of 5-coloring then always 
$ker(\psi)  \neq \{0\}.$}.

We can construct analogous examples if
$gcd(p-1,n) > 2$ \cite{D-J-P}.

\section{Unknotting number from a skein module deformation of moves}

We discuss, in this section, how certain substitutions in 
skein module deformation of elementary moves define weighted 
Fox colorings and can be used to find unknotting numbers of 
links\footnote{We follow the Traczyk's idea how to use 
invariants of links satisfying
skein relations to approximate the unknotting numbers of links.
 Traczyk described this for the Jones polynomial $V_L(e^{2\pi i/6})$
 \cite{Tr,P-2} and in this case $|Col_3(L)|=3|V_L(e^{2\pi i/6})|^2$.
 A.Stoimenov implemented the method for
the Brandt-Lickorish-Millett-Ho polynomial \cite{Sto}.}.

\begin{figure}[ht!]
\centerline{\includegraphics[width=.8\hsize]{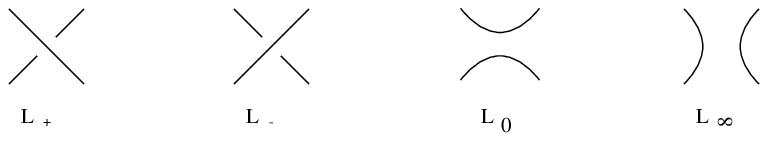}}
\nocolon\caption{}
\end{figure}

We describe here the case of the 2-move whose deformation leads to the 
Kauffman polynomial of links and speculate about the deformation of a 
3-move, cubic polynomial and weighted 7-colorings.

In order to find a relation between the Kauffman polynomial and the
Fox 5-colorings, we study first how the Kauffman polynomial
is  changed under $\pm (2,2)$-moves\footnote{The Kauffman polynomial 
is the Laurent polynomial
 invariant of framed links satisfying the skein relation
$F_{L_+} ( a,x) + F_{L_-}(a,x) =
x(F_{L_0}(a,x) + F_{L_{\infty}} (a,x))$ (Figure  3.1), and the framing relation
$F_{L^{(1)}}(a,x)=aF_L(a,x)$, where $L^{(1)}$ denotes a link obtained
from $L$ by one positive twist of the framing of $L$.}. 
We look for substitutions for $a$ 
and $x$ such that 
\begin{condition}\label{3.1}
$F_L{
\lower .2cm \hbox{\includegraphics[height=.6cm]{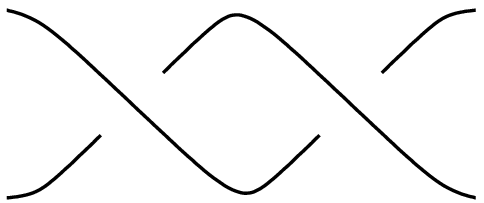}}
} 
= cF_L{
\lower .5cm \hbox{\includegraphics[height=1.2cm]{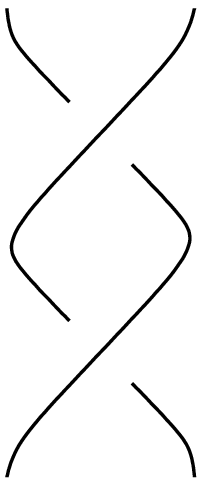}}
}$ 
and 
$F_L{
\lower .2cm \hbox{\includegraphics[height=.6cm]{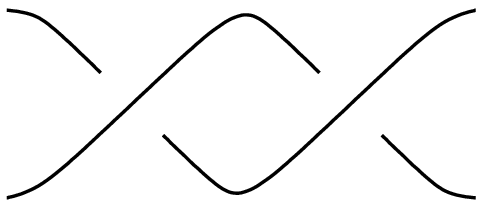}}
} 
= c'F_L{
\lower .5cm \hbox{\includegraphics[height=1.2cm]{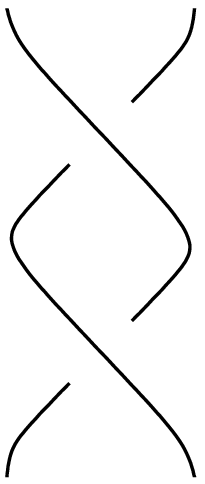}}
}$ for some constants $c$ and $c'$. 
\end{condition}
From the definition we get:

$F_L{
\lower .2cm \hbox{\includegraphics[height=.6cm]{twist2.eps}}
} = xF_{L_+} - F_{L_0} +a^{-1}xF_{L_{\infty}}$,
and $F_L{
\lower .5cm \hbox{\includegraphics[height=1.2cm]{twist-h.eps}}
} = -xF_{L_+} + (a^{-1}x+x^2)F_{L_0} +
(-1+ x^2)F_{L_{\infty}}$.

If we assume Condition 3.1
and compare coefficients of elementary tangles ($L_+,L_0,L_{\infty}$)
we obtain: $c=c'=-a^2=-1$ and $x^2 +ax-1=0$.
These lead to $x= \frac{-a\pm \sqrt{5}}{2}$ (if we let $x=b+b^{-1}$ then
$x^2 +ax-1= b^2 + ab +1 +ab^{-1} + b^{-2} = b^{-2}(\frac{b^5 - a}{b-a})$ 
so $b$ is a 5th primitive root of $1$ for $a=1$, and $10$th primitive root 
of $1$ for $a=-1$). 
Let us assume that $a=1$ and $x=\frac{-1 + \sqrt{5}}{2}=2cos(2\pi/5)$.
Then $F_{T_n}=(\sqrt{5})^{n-1}$ and for links $(2,2)$-equivalent 
to trivial links one has immediately that 
$Col_5(L) = 5(F_L(1,2cos(2\pi/5))^2$. To see this one should compare initial
data (for trivial links) and use  the fact that both sides of the equation
are preserved by $\pm(2,2)$-moves.\  
In fact the equality holds in general.
\begin{lemma}\label{3.2}
For every link $L$ we have $Col_5(L) = 5(F_L(1,2cos(2\pi/5))^2$.
\end{lemma}
\begin{proof}
As mentioned before $Col_5(L)= H_1(M_L^{(2)};Z_5)\oplus Z_5$.
Furthermore Jones proved (\cite{Jon,Ron}), analyzing
the Goeritz matrix of $L$, that
$(F_L(1,2cos(2\pi/5))^2$ is the order of $H_1(M_L^{(2)};Z_5)$.
The diagrammatic proof of Lemma 3.2 was given in \cite{Ja-P}
but due to the untimely death of Francois Jaeger the paper is
not published yet (compare the remark on p. 283 in \cite{P-2}).
\end{proof}

From the fact that every $\pm(2,2)$-move is changing the sign of
$F_L(1,2cos(2\pi/5))$ and from Theorem 3.4 we obtain the following
very elementary but very powerful statement (we formulate it before
Theorem 3.4 because of its very elementary character)\footnote{Our 
method is related to that used by J\,R Rickard, a student of
W\,B\,R Lickorish in his unfinished PhD thesis
where he shows that the knot $8_{16}$ has the unknotting number equal to 2.
Rickard died from cancer before finishing his PhD thesis. 
The short outline
of his work is given in \cite{Lic,Kir}. }:
\begin{corollary}{\rm\cite{P-T-V}}\qua\label{3.3}
\begin{enumerate}
\item[\rm(i)] 
If a knot $K$ can be reduced to the trivial link of two components
 by an even number
of $\pm(2,2)$-moves then the unknotting number $u(K)$ is at least 
two.
\item[\rm(ii)] If a knot $K$ can be reduced to the trivial link of $n$
components by $k$\ $\pm(2,2)$-moves then 
$u(K) \geq n + \frac{(-1)^{n-k} -1}{2}$.
\end{enumerate}
\end{corollary}
From Corollary 3.3 we get that $u(7_4)=2$, $u(8_8)=2$ and
$u(8_{16}) = 2$. \ \ 
Namely, first we check by inspection that
two crossing changes trivialize our knots. In Figure 3.2 we
illustrate how to reduce the knot $7_4$ to $T_2$ by four
$\pm(2,2)$-moves (two 5-moves). In Figure 3.3 the reduction
of the knot $8_8$ by two $\pm(2,2)$-moves is illustrated.
Finally in Figure 3.4 
the reduction of the knot $8_{16}$  by an even number (six)
of $\pm(2,2)$-moves is illustrated.

\begin{figure}[ht!]
\centerline{\includegraphics[width=.7\hsize]{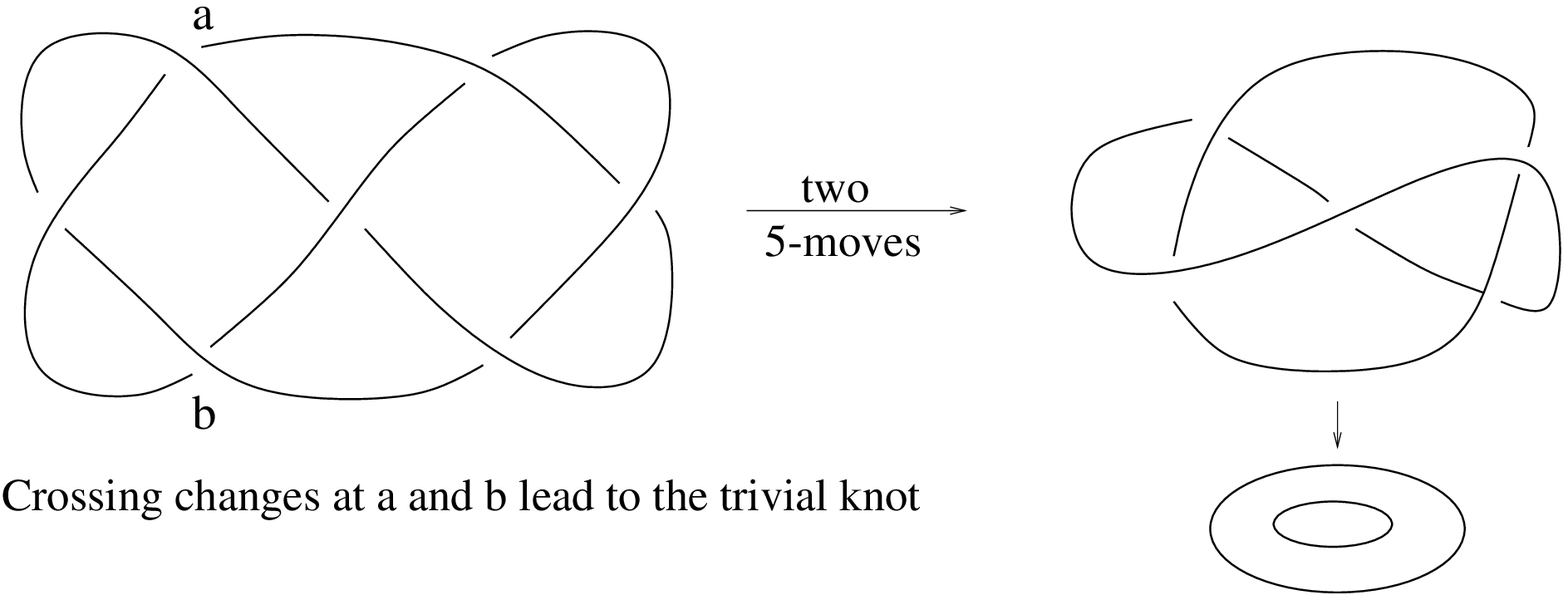}}
\nocolon\caption{}
\end{figure}

\begin{figure}[ht!]
\centerline{\includegraphics[width=.7\hsize]{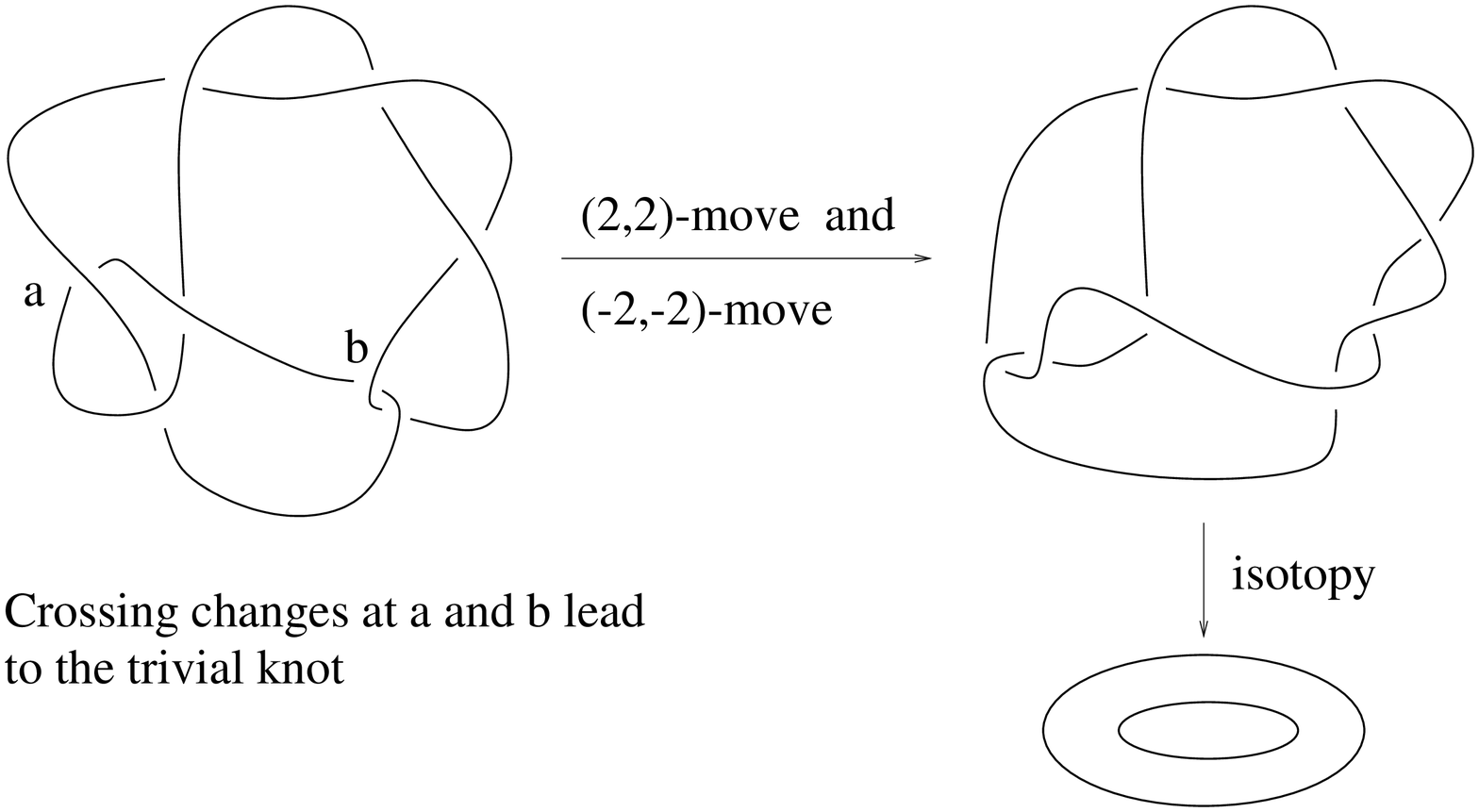}}
\nocolon\caption{}
\end{figure}

\begin{figure}[ht!]
\centerline{\includegraphics[width=.7\hsize]{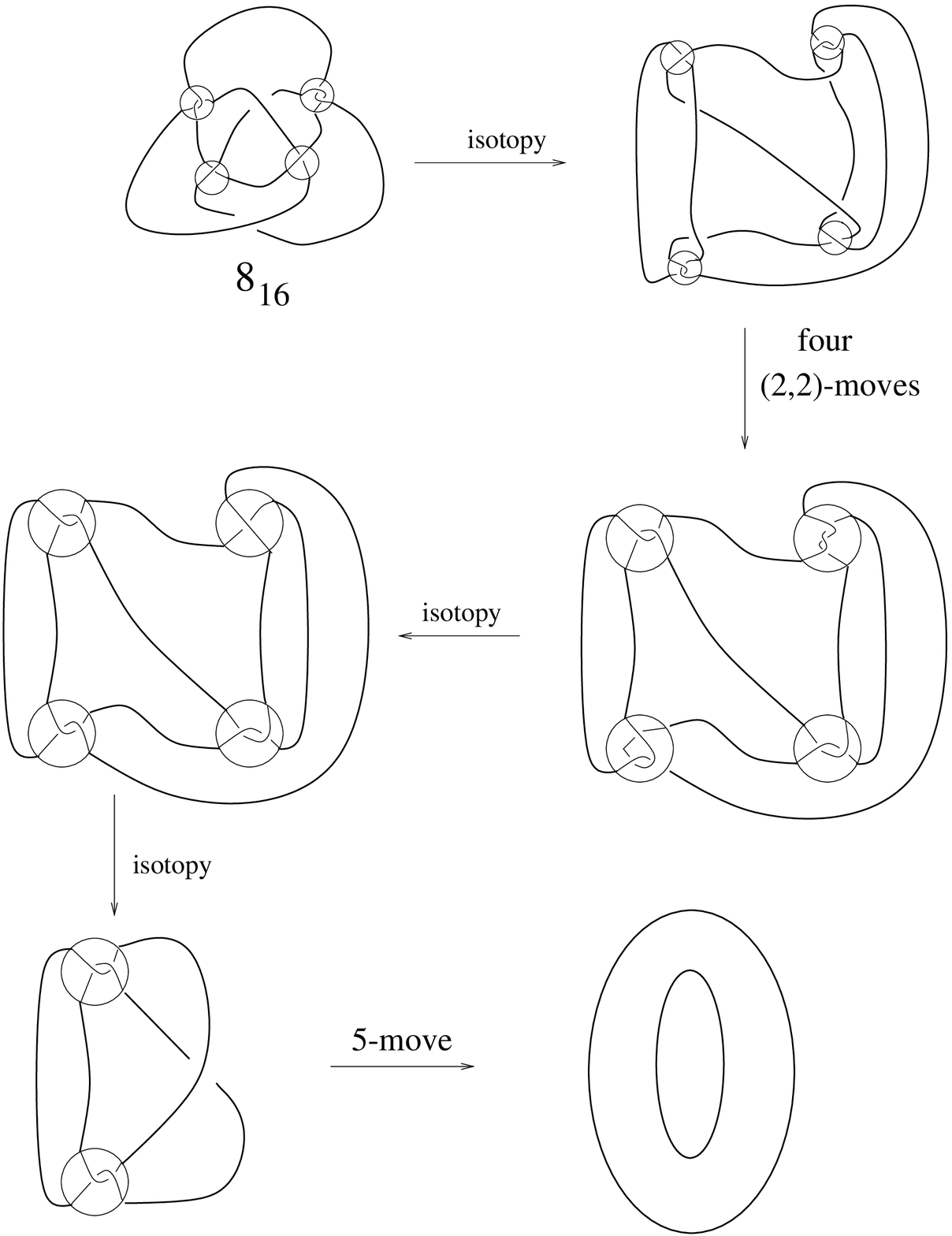}}
\nocolon\caption{}
\end{figure}

For the knot $9_{49}$ one should use Traczyk method 
in full its generality\footnote{The authors of \cite{P-T-V}  
tried, but failed to simplify the knot $9_{49}$ to a trivial link by 
$(2,2)$-moves. Hence, Corollary 3.3 could not be used. The full 
reason for this is explained in Theorem 4.5(3).}.
As we noted before, for $a=1, x=2cos(2\pi/5)$, 
we have $F_L=\pm\sqrt{\frac{1}{5}col_5(L)}=
\epsilon(L)\sqrt{5}^{\lambda (L)}$ where $\epsilon(L)=\pm 1$ and
$\lambda (L)$ is an integer.
\begin{theorem}{\rm\cite{P-T-V,Sto}}\qua\label{3.4}
For a knot $K$, its Kauffman skein quadruplet, $K_+,K_-,$ $K_0,K_{\infty}$ 
and invariants $u(K)$, $\epsilon(K)$ and ${\lambda}(K)$ we have the following
properties of the unknotting number:
\begin{enumerate}
\item[\rm(i)] $u(K) \geq \lambda(K)$.
\item[\rm(ii)] If $u(K_+) = {\lambda}(K_+)$ and $u(K_-)<u(K_+)$
then $u(K_-)={\lambda}(K_-)={\lambda}(K_+)-1$ and
$\epsilon(K_-)=-\epsilon(K_+)=\epsilon(K_0) =\epsilon(K_{\infty})$.
\item[\rm(iii)] If $\epsilon(K) = -(-1)^{{\lambda}(K)}$ then
   $u(K) > \lambda(K)$.
\end{enumerate}
\end{theorem}
\proof
Part (i) is a special case of the Wendt theorem, the first nontrivial
result concerning unknotting number \cite{Wen} 
(compare Lemma 2.2(h) of \cite{P-2})\footnote{$\lambda(K) =log_5Col_5(K)-1$ 
and a crossing change can change dimension of $Col_5(K)$ at most by one.
This can be proven quickly by observing that the relations for $L_+=$
\lower .2cm \hbox{\includegraphics[height=.6cm]{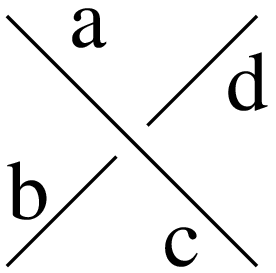}}
are $a=c$ and $a-b+c-d=0$, and the relations for $L_-=$
\lower .2cm \hbox{\includegraphics[height=.6cm]{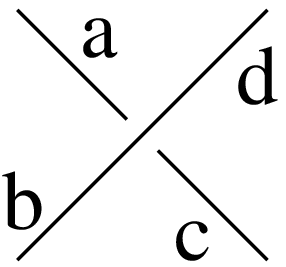}}
are $b=d$ and $a-b+c-d=0$. Thus the linear spaces $Col_5(L_+)$ and 
$Col_5(L_-)$ differ by at most one equation, so their dimensions differ
by at most one.}. 
Assumptions of part (ii) guarantee that ${\lambda}(K_-)={\lambda}(K_+)-1$
as ${\lambda}(K)$ can  be changed at most by 1 when $K$ is modified by
a crossing change. Thus $u(K_-)={\lambda}(K_-)$. Now we have two cases
to consider:
\begin{enumerate}
\item[(a)] $\epsilon(K_+) = \epsilon(K_-)$. \ 
We put $\epsilon(K_+) = \epsilon(K_-) =1$ (the case of $-1$ being analogous).
Then we have $$\sqrt{5}^{{\lambda}(K_+)} + \sqrt{5}^{{\lambda}(K_-)}=
(\frac{-1 + \sqrt{5}}{2})(\epsilon(K_0)\sqrt{5}^{{\lambda}(K_0)} + 
\epsilon(K_{\infty})\sqrt{5}^{{\lambda}(K_{\infty})})$$ 
Therefore
$(\sqrt{5} + 1)\sqrt{5}^{{\lambda}(K_-)}= (\frac{-1 + \sqrt{5}}{2})
(\epsilon(K_0)\sqrt{5}^{{\lambda}(K_0)} + 
\epsilon(K_{\infty})\sqrt{5}^{{\lambda}(K_{\infty})})$
Thus $(3+\sqrt{5})\sqrt{5}^{{\lambda}(K_-)}=
\epsilon(K_0)\sqrt{5}^{{\lambda}(K_0)} + 
\epsilon(K_{\infty})\sqrt{5}^{{\lambda}(K_{\infty})}$
which is impossible (just compare coefficients of $1$ and $\sqrt{5}$ in the
formula). 
\item[(b)] $\epsilon(K_+) = -\epsilon(K_-)$. This case leads to part (iii).
Furthermore, the formula we obtain is as follows:
$$\epsilon(K_+)(\sqrt{5}-1)\sqrt{5}^{{\lambda}(K_-)} 
= (\frac{-1 + \sqrt{5}}{2})
(\epsilon(K_0)\sqrt{5}^{{\lambda}(K_0)} + 
\epsilon(K_{\infty})\sqrt{5}^{{\lambda}(K_{\infty})})$$
Thus $2\epsilon(K_+)\sqrt{5}^{{\lambda}(K_-)} =
(\epsilon(K_0)\sqrt{5}^{{\lambda}(K_0)} +
\epsilon(K_{\infty})\sqrt{5}^{{\lambda}(K_{\infty})})$
Which holds iff  
${\lambda}(K_-)={\lambda}(K_0)= {\lambda}(K_{\infty})$ and
$-\epsilon(K_-)= \epsilon(K_+)=\epsilon(K_0) =\epsilon(K_{\infty})$.
This completes our proof of Theorem 3.4.\endproof
\end{enumerate}

\begin{corollary}\label{3.5}
 The knot $9_{49}$ has the unknotting number equal to 3.
\end{corollary}
\begin{proof}
We have $F_{9_{49}}(1,\frac{-1 + \sqrt{5}}{2})=-5=-\sqrt{5}^2$. Thus
$\epsilon(9_{49})=-1$ and by Theorem 3.2(iii), $u(9_{49}) >2$.
We can easily check by inspection that the knot $9_{49}$ can 
be unknotted by 3 crossing changes thus $u(9_{49})=3$.
\end{proof} 

It is convenient to reformulate Theorem 3.4 so it applies to links.
Consider, after \cite{Mur}, the metric space of links of $n$ components,
with the distance, $u(L_1,L_2)$, defined to be the minimal number
of crossing changes needed to convert $L_1$ to $L_2$.

\begin{theorem}{\rm\cite{P-T-V}}\label{3.6}
$$u(L_1,L_2) \geq |\lambda(L_2)- \lambda(L_1)| + 
\frac{|\epsilon(L_1)\epsilon(L_2) - (-1)^{\lambda(L_2)- \lambda(L_1)}|}{2}.$$
\end{theorem}

\begin{example}\label{3.7}
The distance between the trefoil knot and the figure eight knot
is equal to two (i.e.\ $u(3_1,4_1)=2$).

We have, of course, that 
$u(3_1,4_1) \leq 2$. Furthermore $F_{3_1}(1,2cos(2\pi/5))=-1$ and
$F_{4_1}(1,2cos(2\pi/5))=-\sqrt 5$. Thus $\epsilon(3_1)= \epsilon(4_1) =-1$, 
$\lambda(4_1) = \lambda(3_1) + 1 =1$. Therefore by Theorem 3.6 on has
$u(3_1,4_1) \geq 1+1=2$ and finally $u(3_1,4_1)=2$. 
\end{example} 

We predict that Fox $7$-colorings should be related with cubic
skein modules of $S^3$. In turn the cubic skein module should be useful to
analyze unknotting number for knots with nontrivial Fox 7-colorings.
 We analyzed with M.Veve cubic skein invariants,
preserved, up to constant, by $\pm (2,3)$ moves. We didn't get
results similar to that for the Kauffman or Jones polynomials. We speculate
that the method does not work because instead of a constant one should
probably use the invariant related to $Z_7$ Witt class of the Goeritz
form of a link. This line of research awaits exploration.

\section{Burnside groups of links}

Five months after my talks in Kyoto, several elementary move conjectures
have been solved, or more precisely, disproved, using Burnside groups
of links. I describe below the joint work with my student Mietek
D{\c a}bkowski. Our tool is a noncommutative analogue of
Fox $n$-colorings which we call the $n$'th Burnside group of a link, $B_L(n)$,
\cite{D-P-1,D-P-2}.

\begin{definition}\label{4.1}
The $n$th Burnside group of a link is the
quotient of the
fundamental group of the double branched cover of $S^3$ with the link as
the branch set divided
by all relations of the form $w^n=1$. Succinctly:\
$B_L(n)=\pi_1(M_L^{(2)})/(w^n)$.
\end{definition}
Notice that for the trivial link of $k$ components, $T_k$, one has
$B_{T_k}(n) = B(k-1,n)$ where $B(k-1,n)$ is the classical Burnside
group of $k-1$ generators and exponent $n$.
 
For practical applications it is very convenient to have a diagrammatic
(and local) description of the groups $\pi_1(M_L^{(2)})$ and $B_L(n)$.
Such a presentation, using the core group idea of Bruck \cite{Bruc}
 and Joice \cite{Joy} was given in \cite{F-R,Wa}\footnote{The presentation
was also found before, in 1976, by Victor Kobelsky a student of O.Viro, 
but this was never published \cite{Vi}.}.

\begin{definition}\label{4.2}
Let $D$ be a diagram of a link $L$. We define 
the associated core group $\Pi^{(2)}_D$ of $D$ by the following presentation:
generators of $\Pi^{(2)}_D$ correspond to arcs of the diagram.
Any crossing $v_s$ yields
the relation $r_s=y_iy_j^{-1}y_iy_k^{-1}$ where $y_i$ corresponds
to the overcrossing
and $y_j,y_k$ correspond to the undercrossings at $v_s$ (see Figure  4.1).
\end{definition}
In the above presentation of $\Pi^{(2)}_L$ one relation can 
be dropped since it is a consequence of others.

\begin{figure}[ht!]
\centerline{\includegraphics[height=1.6cm]{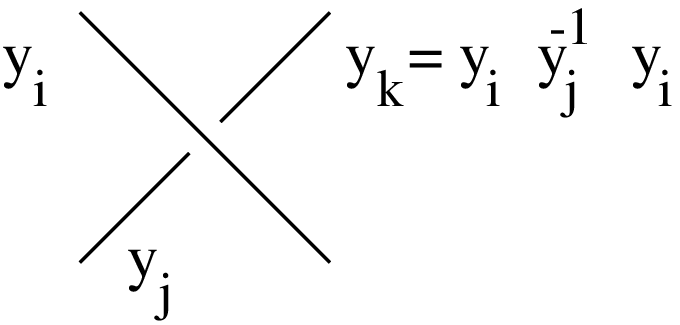}}
\nocolon\caption{}
\end{figure}

The relation to the fundamental group of a double branch cover, mentioned
before, is formulated below (an elementary proof, working 
for links and tangles, using only Wirtinger presentation 
was given in \cite{P-2}).

\begin{theorem}[Wada]\label{4.3}
 $$\Pi^{(2)}_D = \pi_1(M_L^{(2)})\ast Z.$$
Furthermore, if we put $y_i=1$ for any fixed generator, then
$\Pi^{(2)}_D$ reduces to $\pi_1(M_L^{(2)})$. 
\end{theorem}

We introduced Burnside groups in order to analyze elementary moves on links.
 
\begin{theorem}{\rm\cite{D-P-2}}\qua\label{4.4} 
$B_L(n)$ is preserved by rational $\frac{n}{q}$-moves. In 
particular $n$-moves preserve $B_L(n)$.
\end{theorem} 

Theorem 4.4 has been used to find obstructions to several
Conjectures described in Section 1.

\begin{theorem}{\rm\cite{D-P-1,D-P-2}}\qua\label{4.5}
\begin{enumerate}
\item[\rm(1)]
The third Burnside groups of the 2-parallel of the Borromean
rings and of the Chen's link are different than third Burnside groups
of trivial links. In particular we have $|B_{Chen\ link}(3)| = 3^{10}$,
$|B_{T_5}(3)|=3^{14}$, $|B_{2BR}(3)|=3^{21}$ and 
$|B_{T_6}(3)|=3^{25}$.\footnote{It is still an open problem
whether the 2-parallel of the Borromean rings is 3-move equivalent
to the Chen's link with an additional trivial component; both links have 
third Burnside groups of order $3^{21}$.}
\item[\rm(2)]
The fourth Burnside group of the ``half" 2-cabling of the Whitehead 
link, $\cal W$,
is different then the group of the trivial link of 3 component. We have:
$|B_{\cal{W}}(4)|=2^{10} \neq 2^{12}= B(2,4)=B_{T_3}(4)$. 
For the Borromean rings we get $|B_{BR}(4)|=2^5$ 
which gives a simple proof of the Nakanishi
result that the Borromean rings are not 4-move equivalent to a trivial
link.  
\item[\rm(3)] The links $9_{40}$ and $9_{49}$ are not $(2,2)$-move 
equivalent to trivial links\footnote{It is still an open problems
whether $9_{40}$ and $9_{49}$ are $(2,2)$-move equivalent or whether
the Burnside group $B_{9_{40}}(5)$ is isomorphic to $B_{9_{49}}(5)$.}. 
We compare the third terms, $L_3(\ )$, of the graded Lie algebra
associated to the lower central series of groups\footnote{The lower
central series of a group $G$ ($G_1=G, G_2=[G,G],...,
G_n= [G_{n-1},G]$) yields the associated graded Lie ring of the group:
$L= L_1 \oplus L_2 \oplus ... \oplus L_i \oplus ...$ where $L_i= G_i/G_{i+1}$.
The Lie bracket in $L$ corresponds to the group bracket $[g,h]=g^{-1}h^{-1}gh$,\cite{Va}.} 
and show that $L_3(B_{9_{40}}(5)) = L_3(B_{9_{49}}(5)) \neq 
L_3(B_{T_3}(5)) = L_3(B(2,5)) = Z_5 \oplus Z_5$.
\item[\rm(4)] 
For $p$ prime , $p>3$, the closure of the 3-braid $(\sigma_1\sigma_2)^6$ 
cannot be reduced to a trivial link by $\frac{p}{q}$-moves.
We show that the obstruction lies in the third term of the graded Lie algebra
associated to the lower central series of the $p$'th Burnside group. 
\end{enumerate}
\end{theorem}
Notice that Theorem 4.5(1) combined with the fact that 3-algebraic
links are 3-move equivalent to trivial links \cite{P-Ts} implies that
the 2-parallel of the Borromean rings and the Chen's link are not
3-algebraic. The fact that not every link is 3-algebraic is new.
The question whether for any $n$ there is a link which is not
$n$ algebraic still remains open. 
Similarly,
Theorem 4.5(3) combined with Lemma 1.7 implies that
the knots $9_{40}$ and $9_{49}$ are not 2-algebraic\footnote{The
first proof that the  knot $9_{49}$ is not 2-algebraic was
by showing that the
2-fold branched cover of $(S^3,9_{49})$ is a hyperbolic 3-manifold. 
In fact, it is the manifold
I suspected from 1983 to have the smallest volume among
oriented hyperbolic 3-manifolds \cite{I-MPT,Kir,M-V}.}. 

We were unable to use our method in the case when the 
abelianization of the
$n$'th Burnside group (i.e.\ $H_1(M^{2}_L,Z_n)$) is a cyclic group.
In particular Nakanishi 4-move conjecture remains
open as well as 2-component version of the Kawauchi 4-move question 
(Problem 1.4(ii)). In the same vain the question whether the knot $8_{18}$
is $(2,3)$-move equivalent to a trivial link remains open.

\rk{Acknowledgments}
I would like to thank organizers of the Kyoto workshop for their hospitality.

\Addresses\recd

\end{document}